\documentclass{itogi20171}

\usepackage[russian]{babel}
\usepackage[cp1251]{inputenc}

\usepackage{graphicx}

\usepackage{comment}
\usepackage{enumerate} 
\usepackage{xcolor} 

\usepackage[unicode,colorlinks,linkcolor=blue,citecolor=red,bookmarksopen,pdfhighlight=/N]{hyperref}

\newtheorem{thm}{Теорема}

\newtheorem{corollary}{Следствие}
\newtheorem{lemma}{Лемма}

\theoremstyle{definition}
\newtheorem{remark}{Замечание}

\renewcommand{\leq}{\leqslant} 
\renewcommand{\geq}{\geqslant}
\newcommand{\RR}{\mathbb{R}} 
\newcommand{\CC}{\mathbb{C}} 
\newcommand{\NN}{\mathbb{N}} 
\newcommand{\rad}{\mathsf{rad}}
\renewcommand{\Re}{{\rm Re\,}}

\DeclareMathOperator{\conv}{\mathsf{co}}
\DeclareMathOperator{\dd}{d}
\DeclareMathOperator{\intr}{\mathsf{int}}
\DeclareMathOperator{\Hol}{\mathsf{Hol}}
\DeclareMathOperator{\Exp}{\mathsf{Exp}}
\DeclareMathOperator{\spf}{\mathsf{spf}}
\DeclareMathOperator{\clos}{\mathsf{cl}}
\DeclareMathOperator{\area}{\mathsf{area}}
\DeclareMathOperator{\bound}{\mathsf{bd}}

\begin{document}

\title[Полнота  экспоненциальных систем в пространствах функций  \dots и площадь]{Полнота  экспоненциальных систем в пространствах функций на компакте или области  и их площадь}

\author[Хабибуллин Булат Нурмиевич]{Б. Н.  Хабибуллин}
\address{Институт математики с вычислительным центром Уфимского федерального исследовательского центра Российской академии наук}
\email{khabib-bulat@mail.ru}
\thanks{Работа выполнена в рамках государственного задания Министерства науки и высшего образования 
Российской Федерации (код научной темы FMRS-2022-0124).}

\author[Кудашева Елена Геннадьевна]{Е. Г. Кудашева}
\address{Башкирский государственный педагогический университет им. М.~Акмуллы}
\email{lena\_kudasheva@mail.ru}




\UDC{517.538.2, 517.547.22,  514.17, 517.574}

\begin{abstract}
Устанавливаются условия  полноты экспоненциальной системы в пространствах функций, непрерывных  на компакте со связным дополнением  и голоморфных во внутренности этого компакта, а также в пространствах голоморфных функций в ограниченной односвязной области в терминах евклидовой площади выпуклой оболочки этого компакта или области, а также некоторых специальных характеристик или плотностей распределений показателей экспоненциальной системы.
\end{abstract}


\maketitle

\tableofcontents

\section{Введение}

\subsection{Некоторые обозначения, понятия и соглашения}
\textit{Пустое множество\/} обозначаем через  $\emptyset$,  $\mathbb N:=\{1,2, \dots\}$ --- множество всех {\it натуральных чисел,\/} 
 $\NN_0:=\{0\}\bigcup \NN=\{0,1, 2, \dots\}$,  $\overline \NN_0:=\NN_0\bigcup \{+\infty\}$ --- {\it расширение\/}  множества $\NN_0$ со стандартным отношением порядка $\leq$ и точной верхней гранью $+\infty:=\sup \NN_0\notin \NN_0$, для которой  неравенства  $n\leq  +\infty$ выполнены при всех $n\in \overline \NN_0$. 
\textit{Множества всех действительных чисел\/} $\RR$ с таким же отношением порядка $\leq$  
рассматриваем и  как \textit{вещественную ось\/}  в \textit{комплексной плоскости\/} $\CC$ с евклидовой  нормой-модулем $|\cdot|$. Порядковое пополнение множества $\RR$ верхней и нижней гранями $+\infty  :=\sup \RR=\inf\emptyset \notin \RR$ и $-\infty  :=\inf \RR=\sup\emptyset\notin \RR$
определяет  \textit{расширенную\/} вещественную ось $\overline \RR:= \RR\bigcup \{\pm\infty\}$, где, в дополнение к стандартным допустимым операциям, 
полагаем  $0\cdot (\pm\infty)= (\pm\infty) \cdot 0:=0$. 
\textit{Величина} $c\in \overline \RR$ рассматривается и как \textit{функция, тождественно равная\/ $c$,} как правило, на  плоскости $\CC$. Символом  $0$,  наряду с  $0\in \overline \RR$, обозначаем и нулевые функции, меры и т.п.  

\textit{Промежутки с концами  $a\in \overline \RR$ и $b\in \overline \RR$} ---  это множества $[a,b]:=\bigl\{x\in \overline \RR\bigm| a\leq x\leq b\bigr\}$ --- \textit{отрезок\/} в $\overline \RR$,
$(a,b]:=[a,b]\setminus a$, $[a,b):=[a,b]\setminus b$, а  
$(a,b):=[a,b)\setminus a$ и $(a,+\infty]$, $[-\infty, b)$ --- \textit{открытые промежутки\/} в $\overline \RR$, образующие базу открытых множеств при $a<b$.  
Используем и обозначение $\RR^+:=[0, +\infty)$  
 для \textit{положительной полуоси\/} с \textit{расширением\/}    $\overline \RR^+:=[0, +\infty]$.
При этом величина  $x\in \overline \RR$ \textit{положительная\/} при $x\in \overline \RR^+$, 
{\it строго\/} положительная при $0\neq x\in \overline \RR^+$, \textit{отрицательная\/} при $x\in -\overline \RR^+$, 
{\it строго\/} отрицательная при $0\neq x\in -\overline \RR^+$, $x^+:=\sup \{0,x\}\in \overline \RR^+$ --- {\it положительную часть\/} величины   $x\in \overline \RR$, $x^-:=(-x)^+\in \overline \RR^+$ --- её {\it отрицательную часть.\/}

$D(r):=\bigl\{z'\in \CC\bigm| |z|<r\bigr\}$ и $\overline D(r):=\bigl\{z'\in \CC\bigm| |z|\leq r\bigr\}$,
а также $\partial \overline D(r):=\overline D(r)\setminus  D(r)$ 
---  соответственно \textit{открытый и замкнутый круги,}  а также \textit{окружность с центром в нуле  радиуса  $r\in \overline \RR^+$.} Для  $S\subset \CC$ через $\clos S$, $\intr S$, $\bound S$ и  $\conv S$
обозначаем соответственно  \textit{замыкание,\/} \textit{внутренность,\/} \textit{границу\/}  и \textit{выпуклую оболочку\/}  множества $S$ в  $\CC$. 
Таким образом, при $0<r\in \overline \RR^+$ имеют место равенства $\overline D(r)=\clos D(r)$ и $\partial \overline D(r) =\bound D(r)$,  но не при $r=0$, поскольку в этом случае это уже не так, а именно:  $\clos D(0)=\clos \emptyset=\emptyset=\bound D(0)\neq \{0\}=\overline D(0)=\partial \overline D(0)$.

Для {\it расширенной числовой функции\/} $f\colon X\to \overline \RR$ через $f^+\colon x\underset{x\in X}{\longmapsto} \bigl(f(x)\bigr)^+$ обозначаем её {\it положительную часть\/}, а через $f^{-}:=(-f)^+$ --- её {\it отрицательную  часть.\/} Если $f=f^+$,  то \textit{функция $f$ положительная,\/}  
а если  $f=-f^-$, то  {\it отрицательная.\/}  Если $X\subset \overline \RR$ и 
для  любых $x_1,x_2\in X$ из $x_1<x_2$ следует нестрогое неравенство $f(x_1)\leqslant f(x_2)$ (соответственно строгое неравенство $f(x_1)< f(x_2)$), то функция $f$ \textit{возрастающая\/}  (соответственно \textit{строго возрастающая}) на $X$;  функция $f$ \textit{убывающая\/}  (соответственно \textit{строго убывающая}) на $X$, если противоположная функция $-f$ возрастающая  (соответственно строго возрастающая) на $X$.

\subsection{Постановка задачи}
Всюду далее через $Z$ обозначаем \textit{распределение точек\/} на комплексной плоскости\/ $\mathbb C$, среди которых  могут быть повторяющиеся. Распределение точек $Z$ однозначно определяется функцией, действующей из $\mathbb C$ в $\overline \NN_0$ и равной в каждой точке $z\in \mathbb C$ количеству  повторений этой точки $z$  в  $Z$. Для такой функции, которую часто называют \textit{функцией кратности\/}  распределения точек  $Z$ \cite[пп.~0.1.2--0.1.3]{Khsur}, 
или его \textit{дивизором,\/} сохраняем то же  самое обозначение $Z$. Другими словами, $Z(z)$ --- это количество вхождений   точки  $z\in \CC$ в распределение точек $Z$ и  пишем $z\in Z$,  если $Z(z)>0$. Распределение точек  $Z$ можно эквивалентным образом трактовать  и как \textit{распределение масс,\/} или меру, со значениями в  $\overline \NN_0$ с тем же самым обозначением  
\begin{equation}\label{Z}
Z(S):=\sum_{z\in S} Z(z)\in \overline \NN_0\quad\text{\it для любого  $S\subset \CC$.}
\end{equation}
Если \textit{считающая радиальная функция} 
\begin{equation}\label{Zr}
Z^{\rad}(r):=Z\bigl(\overline D(r)\bigr)\overset{\eqref{Z}}{=}
\sum_{z\in \overline D(r)} Z(z)\in \overline \RR^+, \quad r\in \overline \RR^+, 
\end{equation}
для  $Z$ конечна при каждом   $r\in \RR^+$,  т.е. $Z^{\rad}(r)<+\infty$ для всех   $r\in \RR^+$, 
то  $Z$  \textit{локально конечное.} 

\textit{Евклидову площадь\/}  множества $S\subset \CC$ обозначаем через 
\begin{equation*}
\area(S):=\iint_S \dd x\dd y=\iint_S r\dd r\dd \theta, \quad   x+iy=re^{i\theta},\; x,y,\theta\in \RR, \; r\in \RR^+,
\end{equation*} 
если двойные интегралы справа корректно определены, или множество $S$ измеримо по плоской мере Лебега на $\CC$. В частности, это  всегда имеет место для \textit{выпуклых ограниченных\/} $S\subset \CC$.

Система векторов из  топологического векторного пространства \textit{полна\/} в нём,  
если замыкание линейной оболочки этой системы совпадает с этим  пространством. Для распределения точек $Z$ на $\CC$ в данной статье  далее исследуется  полнота лишь  \textit{экспоненциальных систем}
 \begin{equation}\label{Exp}
\Exp^Z:=\Bigl\{ w\underset{w\in \CC}{\longmapsto}w^p\exp (zw)\Bigm| z\in Z, \; Z(z)-1\geq p\in \NN_0\Bigr\}
\end{equation}
\textit{с распределением показателей $Z$,} что, в частности, актуально в спектральной теории операторов.  

Для функции $f$ на  $S\subset \CC$  со значениями в $\CC$ или в $\overline \RR$
полагаем 
\begin{equation}\label{CKn}
\|f\|_{S}:=\sup\Bigl\{\bigl|f(z)\bigr|\Bigm| z\in S \Bigr\},
\end{equation}
а через $C(S)$ обозначаем \textit{пространство непрерывных функций $f\colon S\to \CC$ с\/  $\sup$-нормой} \eqref{CKn}.
Для \textit{открытого\/} подмножества $S\subset \CC$ через $\Hol(S)$ обозначаем \textit{пространство голоморфных функций $f\colon S\to \CC$ с топологией равномерной сходимости на всех компактах $K\subset S$,\/} определяемой $\sup$-полунормами $\|f\|_{K}$. Для \textit{компакта} $S\subset \CC$ с \textit{внутренностью\/} $\intr S$ через $C(S)\bigcap \Hol(\intr S)$  обозначаем \textit{банахово пространство  непрерывных на $S$ и голоморфных на внутренности $\intr S$ функций $f\colon S\to \CC$ с $\sup$-нормой\/} $\|f\|_{S}$.  Таким образом,  если в последнем случае $\intr S=\emptyset$ --- пустое множество, то $C(S)\bigcap \Hol(\intr S)$ --- это банахово пространство $C(S)$ непрерывных на $S$ функций со значениями в $\CC$. 
Книга \cite[п.~3.2]{Khsur}  содержит детальный обзор по вопросам полноты экспоненциальных систем $\Exp^Z$ по состоянию  вплоть до 2012 г.  в разнообразных функциональных пространствах ---  в значительной мере именно для пространств  $\Hol(S)$ или  $C(S)\bigcap \Hol(\intr S)$ функций соответственно на области $S\subset \CC$ или компакте $S\subset \CC$. 
\textit{Основная задача\/} --- получить условия  полноты экспоненциальной системы   $\Exp^Z$ из \eqref{Exp} в 
функциональных пространствах $\Hol(S)$ или  $C(S)\bigcap \Hol(\intr S)$, когда соответственно для ограниченной области $S$ или компакта $S$ априори известна лишь  \textit{евклидова площадь\/} $\area(\conv S)$ его \textit{выпуклой оболочки\/} $\conv S$ или площадь $\area(S)$ в случае \textit{выпуклой\/} соответственно области или компакта  $S\subset \CC$.
Естественно требовать, чтобы эти условия выражались через соотношения между какими-либо характеристиками распределения точек-показателей $Z$ и площадью $\area(\conv S)$ и были точны. В данной работе мы не останавливаемся на подтверждении  точности наших результатов, хотя это так и есть, например, для \textit{любых выпуклых\/} $S$. Требуемые для этого примеры основаны на  построении довольно тонких примеров целых функций экспоненциального типа и очень регулярного роста, особенно для компактов $S$.  Построение таких примеров предполагается обсудить в другой работе. 

\subsection{Основные результаты}

В теории целых функций одной комплексной переменной  \cite{Levin56}, \cite{Levin96}, \cite[гл.~I, \S~1, п.~4]{Leo}
а также  в некоторых других вопросах, связанных  с геометрией на плоскости \cite[отдел третий, гл.~3, \S~1]{PolSeg78}, \cite[гл. 1, \S~2]{Sant},   \textit{опорную функцию\/} подмножества $S\subset \CC$ чаще всего  определяли как  функцию 
\begin{equation}\label{spf}
{\spf}_S\colon \theta\underset{\theta \in \RR}{\longmapsto} \sup_{s\in S}\Re se^{-i\theta}\in \overline \RR.
\end{equation} 
По построению  \eqref{spf}   $2\pi$-периодические на $\RR$  опорные функции множества $S\subset \CC$,
его  замыкания $\clos S$ в $\CC$, выпуклых оболочек  $\conv S$, $\conv \clos S$ 
и  замыкания $\clos \conv  S$ совпадают. Сдвиг компакта или ограниченной области $S$ в $\CC$ не влияет на полноту экспоненциальной системы $\Exp^Z$ соответственно в пространстве $C(S)\bigcap \Hol(\intr S)$ или $\Hol(S)$. Поэтому, не умаляя общности, всюду далее нам удобно считать, что, после сдвига ограниченного  $S$ и сохранении за ним того же обозначения $S$, нулевая точка принадлежит замыканию выпуклой оболочки множества $S$, т.е. выполнено условие  
\begin{equation}\label{Kh0S}
0\in \clos \conv S.
\end{equation} 
Если $S\subset \CC$ --- \textit{компакт,\/} то условие  \eqref{Kh0S} \textit{эквивалентно\/} условию $0\in \conv S$.  Кроме того, по определению \eqref{spf}
\textit{для произвольного $S\subset \CC$ условие \eqref{Kh0S} эквивалентно  положительности опорной функции}
\begin{equation}\label{spfgeq0}
\spf_S(\theta)\overset{\eqref{Kh0S}}{\underset{\theta\in \RR}{\geq}} 0. 
\end{equation}

При  трактовке распределения точек $Z$ как распределения масс   \eqref{Z}
 для любой положительной функции $f$ на $S$ можно корректно  определить сумму 
\begin{equation}\label{sumZ}
\sum_{\stackrel{z\in Z}{z\in S}}f(z):=\int_Sf\dd Z \in \overline \RR^+.
\end{equation}

Для точки  $z\in \CC$ через  $\arg z\subset \RR$ обозначаем  множество значений всех её \textit{угловых аргументов.\/} 
Для \textit{$2\pi$-периодической функции\/}  на $\RR$ со значениями в $\RR$ однозначно определены значения этой функции на $\arg z$ для любых $z\in \CC\setminus \{0\}$. При ограниченном $S\subset \CC$ \textit{считающую   радиальную   функцию  для распределения точек  $Z$  по аргументам
относительно $S$\/} определяем как 
\begin{equation}\label{Zradk}
Z^{\rad}_S(r)\underset{r\in \RR^+}{\overset{\eqref{sumZ}}{:=}}
Z(0)\|\spf_S\|_{\RR}
+\sum_{\stackrel{z\in Z}{0<|z|\leq r}}\spf_S(\arg z) \in \overline \RR^+, 
\quad \|\spf_S\|_{\RR}{\overset{\eqref{CKn}}{:=}}\sup_{\theta\in \RR}\,\bigl|\spf_S(\theta)\bigr|.
\end{equation}
При условии \eqref{Kh0S}--\eqref{spfgeq0}  функция $Z^{\rad}_S$ \textit{положительная, возрастающая  и непрерывная  справа на $\RR^+$.} 

Для   $S:=D(r)$ или $S:=\overline D(r)$  их опорные  функции тождественно равны 
\begin{equation}\label{S=D}
\spf_{D(r)}(\theta)\underset{\theta\in \RR}{\equiv}\spf_{\overline D(r)}(\theta)\underset{\theta\in \RR}{\equiv} r\quad\text{и} \quad  
Z^{\rad}_{D(1)}=Z^{\rad}_{\overline D(1)}\overset{\eqref{Zr}}{=}Z^{\rad} 
\end{equation} 
---  считающая радиальная функция $Z^{\rad}$ из \eqref{Zr}. 
В отличие от последней считающая  радиальная  функция  для $Z$ по аргументам
относительно ограниченного $S\subset \CC$ учитывает распределение точек из $Z$ не только по радиусу, но и по аргументам. 

\textit{Комплексно сопряжённое\/} к  $z=re^{i\theta}\in \CC$ с $r\in \RR^+$ и  $\theta \in \RR$ число обозначаем 
через $\bar{z}:=re^{-i\theta}$, а для подмножества $S\subset \CC$   {\it сопряжённое\/} подмножество, зеркально симметричное  $S$ относительно  вещественной оси $\RR$, обозначаем как $\bar{S}:=\bigl\{\bar{z}\bigm| z\in S\bigr\}$ с соответственно   опорной  функцией $\spf_{\bar S}$.

Следующая теорема  ---  частный случай результата,  анонсированного  в \cite[основная теорема]{Kha23V}.

\begin{thm}\label{th1} 
Пусть    $Z$ --- распределение точек в $\CC$. Если для компакта  $S\subset \CC$   со связным дополнением $\CC\setminus S$ при оговорённом в \eqref{Kh0S}--\eqref{spfgeq0}  условии $0\overset{\eqref{Kh0S}}{\in} \clos \conv S=\conv S$  для  некоторого строго положительного $r_0\in \RR^+$ выполнено равенство 
\begin{equation}\label{sZar} 
\sup_{r_0\leqslant r< R<+\infty}\biggl(\int_{r}^{R}\frac{Z_{\bar S}^{\rad}(t)}{t^2}\dd t
- \frac{\area(\conv S)}{\pi}\ln\frac{R}{r}\biggr)=+\infty, 
\end{equation}
то  система $\Exp^Z$ из \eqref{Exp} полна в пространстве $C(S)\bigcap \Hol (\intr S)$. 
\end{thm}

\begin{corollary}\label{cor1}
Если для  произвольных  распределения точек  $Z$ в $\CC$ и  ограниченной односвязной области  $S\subset \CC$   при условии \eqref{Kh0S}--\eqref{spfgeq0} 
выполнено неравенство 
\begin{equation}\label{sZarL}
\limsup_{1<a\to +\infty}\frac{1}{\ln a}\limsup_{0<r\to +\infty} \int_{r}^{ar}\frac{Z_{\spf_{\bar S}}^{\rad}(t)}{t^2}\dd t \geq \frac{1}{\pi} \area(\conv  S), 
\end{equation}
то  система $\Exp^Z$ из \eqref{Exp} полна в пространстве $\Hol (S)$. 
\end{corollary}

\begin{remark} Величину в левой части неравенства \eqref{sZarL}
 по аналогии с подобными плотностями   из  \cite{MR}, \cite[гл. 22]{RC}, \cite[гл. 3]{Khsur}, \cite{KKh00}, \cite{SalKha20}, \cite{SalKha21} можно назвать \textit{верхней логарифмической блок-плотностью для $Z$ относительно евклидовой площади выпуклой оболочки множества $\bar S$.}
\end{remark}

\section{Доказательства результатов}

\begin{proof}[Доказательство теоремы\/ {\rm \ref{th1}}.] Предположим противное, а именно: в условиях теоремы \ref{th1} система $\Exp^Z$ не полна в пространстве $C(S)\bigcap \Hol(\int S)$. Тогда по теореме Рисса о представлении линейных непрерывных функционалов на пространстве $C(S)$ вкупе с теоремой Хана\,--\,Банаха о продолжении линейных непрерывных функционалов с сохранением нормы --- в данном случае с замкнутого подпространства 
$C(S)\bigcap \Hol(\intr S)$ в $C(S)$,  а также из известных следствий из неё,  существует борелевская комплекснозначная мера $\mu\neq 0$ с носителем в $S$, аннулирующая экспоненциальную систему $\Exp^{Z}$, но не аннулирующая хотя бы одну  функцию из $C(S)\bigcap \Hol(\int S)$.  Последнее означает, что заданная преобразованием Фурье\,--\,Лапласа целая функция 
\begin{equation}\label{F}
F\colon z\underset{z\in \CC}{\longmapsto} \int_{S} e^{zs}\dd \mu(s)
\end{equation} 
\textit{обращается в нуль на\/} $Z$ с учётом кратности, а именно: кратность корня функции $F$ в  каждой точке $z\in \CC$  не превышает $Z(z)$. 
В силу  связности дополнения $\CC\setminus S$ целая функция $F$ из \eqref{F} ненулевая. Действительно, если $F=0$, то по \eqref{F} мера $\mu$ аннулирует экспоненциальную систему $\bigl\{e^{zs}\bigm| z\in \CC\bigr\}$, замыкание линейной оболочки которой содержит все многочлены \cite[терема 1 (критерий)]{Gro03}, \cite[гл. 1, п.~1.1.1, пример 1.1.1]{Khsur}. Следовательно, мера $\mu$ аннулирует все многочлены. Но по теореме-критерию Мергеляна при условии связности $\CC\setminus S$  множество всех многочленов плотно в   $C(S)\bigcap \Hol(\int S)$. Тогда мера $\mu$ аннулирует все функции из  $C(S)\bigcap \Hol(\int S)$, что не согласуется с выбором меры $\mu$.  Таким образом, далее $F\neq 0$. 

В обозначении $|\mu|$ для \textit{полной вариации меры\/} $\mu$ и записи  $z:=re^{i\theta}$ в полярной форме с $r\in \RR^+$ и $\theta\in \RR$ для целой функции $F\neq 0$
из \eqref{F} имеет место оценка сверху 
\begin{equation*}
\bigl|F(re^{i\theta})\bigr|\underset{re^{i\theta}\in \CC}{\leq} \int_S\bigl|\exp({re^{i\theta}s})\bigr|\dd |\mu|(s)
\underset{re^{i\theta}\in \CC}{\leq} \sup_{s\in S}\bigl|\exp (re^{i\theta}s )\bigr| |\mu| (S)
\underset{re^{i\theta}\in \CC}{=}\exp \Bigl( r \sup_{s\in S}\Re s e^{i\theta}\Bigr) |\mu| (S),
\end{equation*} 
где по определению опорной функции \eqref{spf}
\begin{equation*}
\sup_{s\in S}\Re s e^{i\theta}\underset{\theta\in \RR}{=}\sup_{s\in S}\Re s e^{-i(-\theta)}\overset{\eqref{spf}}{\underset{\theta\in \RR}{=}}
\spf_{\bar S}(\theta)
\end{equation*}
--- значения опорной функции сопряжённого компакта $\bar S$ в точках $\theta\in \RR$.  
Следовательно, эта  оценка после логарифмирования может быть продолжена  как 
\begin{equation}\label{Fspf}
 \ln\bigl|F(re^{i\theta})\bigr|\underset{re^{i\theta}\in \CC}{\leq}\spf_{\bar S}(\theta)r+\ln |\mu| (S), \quad  |\mu| (S)\neq 0.
\end{equation}
Поскольку целая функция $F$ \textit{ненулевая,}  можем рассмотреть 
субгармоническую функцию 
\begin{equation}\label{Fspfu}
u:=\ln|F|-\ln |\mu| (S)\neq -\infty
\end{equation} 
 с распределением масс Рисса, или мер Рисса  \cite{HK}, \cite{Hor94}, \cite{Rans}
\begin{equation}\label{ulnFm}
\varDelta_u:=\frac{1}{2\pi}{\bigtriangleup}u=\frac{1}{2\pi}{\bigtriangleup}\ln|F|\geq 0, \quad\text{где $\bigtriangleup$ --- \textit{оператор Лапласа,}} 
\end{equation} 
действующий на субгармоническую  функцию $u$ как обобщённую функцию на  пространстве основных финитных функций на $\CC$. 
В частности, ввиду обращения в нуль целой функции $F$ на $Z$ и известного вида   \cite[теорема 3.7.8]{Rans} распределения масс Рисса субгармонической функции $\ln|F|$ при рассмотрении распределения точек $Z$ как распределения масс в смысле \eqref{Z} имеет место неравенство  
$\frac{1}{2\pi}{\bigtriangleup}\ln|F|\geq Z$ и, как следствие, приходим к неравенству $\varDelta_u\overset{\eqref{ulnFm}}{\geq} Z$ для распределений масс 
$\varDelta_u$ и $Z$ на $\CC$. Это неравенство в силу положительности  $\spf_{\bar S}\geq 0$ при условии \eqref{Kh0S} 
по  условию-равенству \eqref{sZar} показывает, что 
\begin{equation}\label{supua}
\sup_{r_0\leq r<R<+\infty}\biggl(\int_r^R\frac{(\varDelta_u)_{\bar S}^{\rad}(t)}{t^2}\dd t- \frac{\area(\bar S)}{\pi}\ln\frac{R}{r}\biggr)=+\infty, 
\end{equation}
где в порядке переноса  определения \eqref{Zradk} с распределений точек на распределения масс положено 
\begin{equation}\label{Zradkdu}
(\varDelta_u)^{\rad}_{\bar S}(t)\underset{r\in \RR^+}{\overset{\eqref{Zradk}}{:=}}
(\varDelta_u)\bigl(\overline D(r_0)\bigr)\|\spf_{\bar S}\|_{\RR}
+\int_{r_0<|z|\leq t}\spf_S(\arg z)\dd (\varDelta_u)(z)\quad\text{при $ t\in [r_0,+\infty)$} 
\end{equation}
--- \textit{считающая   радиальная   функция  для распределения масс  $\varDelta$  по аргументам
относительно $\bar S$ вне $\overline D(r_0)$,\/} которая при \eqref{Kh0S}--\eqref{spfgeq0}  \textit{положительная, возрастающая  и непрерывная  справа на $[r_0,+\infty)$.} 
Напомним, что для \textit{субгармонической на $\CC$}  и непрерывной  функции \cite{GriMal15}
\begin{equation}\label{Ms}
M\colon re^{i\theta}\underset{re^{i\theta}\in \CC}{\longmapsto} \spf_{\bar S}(\theta)r 
\end{equation}
её  распределение масс Рисса определяется   как произведение мер \cite[п. 3.3.1]{Khsur} через её плотность 
\begin{equation}\label{rl}
\dd \varDelta_M(re^{i\theta})=\frac{1}{2\pi}\dd r \otimes \dd l_{\conv \bar S}(\theta)
\end{equation}
в полярных координатах, где $l_{\conv \bar S}(\theta)$ --- длина дуги границы $\bound \conv \bar S$, отсчитываемой при движении по границе <<против часовой стрелки>> от последней точки опоры   опорной к компакту $\clos \conv \bar S$ прямой, ортогональной положительной полуоси $\RR^+$, до последней точки опоры   опорной к компакту $\clos \conv \bar S$ прямой, ортогональной направлению радиус-вектора  точки $e^{i\theta}$ \cite{BF}, \cite{Sant}, \cite[п. 3.3.1]{Khsur}. В частности, вычисление площади выпуклого компакта $\conv \bar S\ni 0$ путём аппроксимации  его выпуклыми описанными многоугольниками, площади которых вычисляются через сумму площадей внутренних треугольников  с центрами  в нулей как половины произведений длин апофем на длины соответствующих сторон, дают равенство для площади 
\begin{equation}\label{areal}
\area(\bar S)=\frac{1}{2}\int_0^{2\pi}\spf_{\bar S}(\theta)\dd l_{\conv \bar S}(\theta).
\end{equation}
Отсюда для вычитаемого произведения с $\area(\bar S)$ из \eqref{supua} при всех $r_0\leq r<R<+\infty$ имеем
\begin{equation}\label{supua1}
 \frac{\area(\bar S)}{\pi}\ln\frac{R}{r}=\frac{1}{2}\int_0^{2\pi}\spf_{\bar S}(\theta)\dd l_{\conv S}(\theta) \int_r^R\frac{1}{t}\dd t  
\overset{\eqref{rl}}{=}\int_{r<|te^{i\theta}|\leq R} \frac{\spf_{\bar S}(\theta)}{r}\dd \varDelta_M(re^{i\theta}), 
\end{equation}
а для  интеграла из \eqref{supua} при всех $r_0\leq r<R<+\infty$ интегрирование  по частям даёт  равенство 
\begin{equation}\label{dui}
\int_r^R\frac{(\varDelta_u)_{\bar S}^{\rad}(t)}{t^2}\dd t
=\frac{(\varDelta_u)_{\bar S}^{\rad}(R)}{R}-\frac{(\varDelta_u)_{\bar S}^{\rad}(r)}{r}
+\int_r^R\frac{1}{t}\dd (\varDelta_u)_{\bar S}^{\rad}(t).
\end{equation}
При этом  в силу \eqref{Fspf} и \eqref{Fspfu} имеют место отграничения 
\begin{equation}\label{uMs}
u(re^{i\theta})\leq \spf_{\bar S}(\theta)r\overset{\eqref{Ms}}{=:}M(re^{i\theta})\quad\text{при всех $re^{i\theta}\in \CC$,}
\end{equation}
откуда $u$ --- \textit{субгармоническая функция  конечного типа при порядке\/} $1$ \cite[гл. 4]{HK}, для которой 
\begin{equation}\label{Du}
\limsup_{z\to \infty}\frac{u(z)}{|z|}<+\infty, \quad \limsup_{r\to +\infty} \frac{\varDelta_u\bigl(\overline D(r)\bigr)}{r}<+\infty. 
\end{equation}
В частности, из последнего предельного соотношения ввиду ограниченности $\spf_{\bar S}$  на $\RR$ получаем 
\begin{equation*}
\sup_{r\in [r_0,+\infty)}\frac{(\varDelta_u)_{\bar S}^{\rad}(r)}{r}\overset{\eqref{Zradkdu}}{\leq} 
\sup_{r\in [r_0,+\infty)} \|\spf_{\bar S}\|_{\RR}\frac{\varDelta_u\bigl(\overline D(r)\bigr)}{r}<+\infty, 
\end{equation*}
откуда для первых двух слагаемых в правой части \eqref{dui} имеем 
\begin{equation*}
\sup_{r_0< r\leq R}\biggl|\frac{(\varDelta_u)_{\bar S}^{\rad}(R)}{R}-\frac{(\varDelta_u)_{\bar S}^{\rad}(r)}{r}\biggr|
\leq 2\sup_{r\in [r_0,+\infty)}\frac{(\varDelta_u)_{\bar S}^{\rad}(r)}{r}<+\infty, 
\end{equation*}
а  последний интеграл в \eqref{dui} согласно \eqref{Zradkdu} можем записать как 
\begin{equation*}
\int_r^R\frac{1}{t}\dd (\varDelta_u)_{\bar S}^{\rad}(t)=
\int_{r<|te^{i\theta}|\leq R}\frac{\spf_{\bar S}(\theta)}{t}\dd \varDelta_u(te^{i\theta}).
\end{equation*}
Отсюда согласно \eqref{supua} и \eqref{supua1} получаем 
\begin{equation}\label{suMd}
\sup_{r_0<r<R<+\infty}\int_{r<|te^{i\theta}|\leq R}\frac{\spf_{\bar S}(\theta)}{t}
\dd \bigl(\varDelta_u-\varDelta_M\bigr)(te^{i\theta})=+\infty.
\end{equation}
Итоговая наша задача --- получить противоречие между этим равенством и ограничением \eqref{uMs}, показав,что из 
\eqref{uMs} следует конечность левой части \eqref{suMd}. 
Для этого рассмотрим функцию 
\begin{multline}\label{VR}
V_R\colon te^{i\theta}\underset{te^{i\theta}\in \CC\setminus \{0\}}{\longmapsto} \spf_{\bar S}(\theta)\Bigl(\frac{1}{t}-\frac{t}{R^2}\Bigr)
\overset{\eqref{spf}}{=}\sup_{s\in S}\bigl(\Re se^{i\theta}\bigr) \Bigl(\frac{1}{t}-\frac{t}{R^2}\Bigr)\\
\underset{z:=te^{i\theta}\neq 0}{=}\; \sup_{s\in S}\Re  \Bigl(\frac{s}{\bar z}-\frac{sz}{R^2}\Bigr)\underset{z\in \CC\setminus\{0\}}=V_R(z),
\end{multline}
которая по построению \textit{положительна на\/} $\overline D(R)$ ввиду условия \eqref{Kh0S}--\eqref{spfgeq0}, \textit{обращается в нуль на окружности\/} $\partial \overline D(R)$, \textit{непрерывна\/} ввиду непрерывности опорных функций ограниченных множеств.  Кроме того, согласно последнему  равенству
в  \eqref{VR}, функция $V_R$ представляет собой  точную верхнюю грань локально ограниченного сверху семейства \textit{гармонических\/} на $\CC\setminus \{0\}$ функций 
\begin{equation*}
\biggl\{\Re  \Bigl(\frac{s}{\bar z}-\frac{sz}{R^2}\Bigr)\biggr\}_{s\in S}. 
\end{equation*}
Отсюда сразу следует, что \textit{функция\/ $V_R$ субгармоническая на\/ $\CC\setminus \{0\}$.}

При этом выпуклый компакт $\conv \bar S\subset \CC$ можно представить как пересечение последовательности выпуклых компактов $K_n\underset{n\in \NN}{\supset} K_{n+1}$, вложенных друг в друга,   для которых их опорные функции $k_n:=\spf_{K_n}$ дважды непрерывно дифференцируемы. По построению убывающая последовательность положительных опорных функций $k_n$  стремится к опорной функции $\spf_{\bar S}$ и  
функции 
\begin{equation}\label{vn}
v_n\colon te^{i\theta}
\underset{te^{i\theta}\in \CC\setminus \{0\}}{\longmapsto} k_n(\theta)\Bigl(\frac{1}{t}-\frac{t}{R^2}\Bigr),
\end{equation}
согласно обоснованному выше,  {\it положительны  на\/}   $\overline D(R)\setminus \{0\}$, а также  \textit{субгармонические\/} и имеют {\it непрерывные частные производные  до второго порядка\/} включительно на $\CC\setminus \{0\}$.

Далее на потребуется следующее объединение двух утверждений из \cite{Kha91_1}, которые могут быть выведены  и  по общим интегральным  формулам  из \cite[теорема 2]{Men22}.

\begin{lemma}[{\rm \cite[леммы 2.2--2.3]{Kha91_1}}]\label{lem2_2} 
Пусть $0<r<R<+\infty$ и  функция $V$ положительна на замкнутом кольце $\overline D(R)\setminus D(r)$,
 субгармоническая в его  внутренности $D(R)\setminus \overline D(r)$, 
тождественно равна  нулю на окружности $\partial \overline D(R)$ и совпадает с  сужением  
на   $\overline D(R)\setminus D(r)$
некоторой дважды непрерывно дифференцируемой  в окрестности  кольца
 $\overline D(R)\setminus D(r)$ функции. 
Используя инверсию  функции $V$ относительно окружности $\partial \overline D(r)$, построим  положительную 
на\/ $\CC$ функцию 
\begin{equation}\label{V*}
V^*(z):=
\begin{cases}
V(z)&\text{при $r< |z|\leqslant R$},\\
V(r^2/\Bar z)&\text{при $r^2/R< |z|\leqslant  r$},\\
0&\text{при  $|z|\leqslant  r^2/R$ и $|z|>R$},
\end{cases} \quad z\in \CC.
\end{equation}
Тогда для любой пары  субгармонических на окрестности круга $\overline D(R)$ функций $u\neq -\infty$ и $M$ с распределениями масс Рисса соответственно $\varDelta_u$ и $\varDelta_M$
из неравенства   $u\leq M$ на этой окрестности следует неравенство 
\begin{equation}\label{VuMn}
\int_{\CC} V^*\dd \varDelta_u\leqslant 
\int_{\CC} V^*\dd \varDelta_M+ \frac{r}{\pi}\int_0^{2\pi} \bigl(u(re^{i\theta})-M(re^{i\theta})\bigr)
\frac{\partial V}{\partial \vec n_{\operatorname{out}}}(re^{i\theta}) \dd \theta, 
\end{equation}  
где  
$\dfrac{\partial}{\partial \vec n_{\operatorname{out}}}$  --- оператор дифференцирования по внешней нормали  
к кольцу $D(R)\setminus \overline D(r)$ на  $\partial \overline D(r)$.
\end{lemma}

Интегральное среднее  функции $g\colon \partial \overline D(r)\to \overline \RR$ по окружности $\partial \overline D(r)$ обозначаем как 
\begin{equation}\label{Cu}
g^{\circ}(r):=\frac{1}{2\pi}\int_0^{2\pi}g(re^{i\theta})\dd \theta.
\end{equation}
Следующая лемма --- предельная форма предшествующей леммы \ref{lem2_2}. 
\begin{lemma}\label{lem2}
Пусть в убывающей последовательности функций $v_n\colon \overline D(R)\setminus D(r)\to \RR$ 
каждая из них  обладает теми  же свойствами, что и функция $V$ в предыдущей лемме, а также 
модули производных по  радиусу от них равномерно по $n$ ограничены сверху во всех точках на окружности $\partial \overline D(R)$ некоторым  числом $N_r\in \RR^+$. Обозначим теперь через $V$ уже предельную функцию для последовательности $v_n$. Тогда для субгармонических на окрестности замкнутого круга $\overline D(R)$ функций $u\leq M$, где $u\neq -\infty$,    имеет место неравенство 
\begin{equation}\label{VuMnv}
\int_{r<|z|\leq R} V\dd \varDelta_u\leqslant \int_{r<|z|\leq R} V\dd \varDelta_M+ \varDelta_M\bigl(\overline D(r)\bigr)\sup_{r\leq |z|\leq R} V(z)  + N_r\frac{r}{\pi}\bigl(|u|^{\circ}(r)+|M|^{\circ}(r)\bigr).
\end{equation}  
\end{lemma}
\begin{proof}[Доказательство леммы\/ {\rm \ref{lem2}}.] 
Производная по внешней нормали 
$\frac{\partial}{\partial \vec n_{\operatorname{out}}}$ 
к кольцу $D(R)\setminus \overline D(r)$ на  $\partial \overline D(r)$ --- это, с точностью до знака, производная по радиусу  на  $\partial \overline D(r)$.  
Поэтому в силу положительности функции $V^*=\lim\limits_{n\to +\infty} v_n^*$, известной теореме о монотонном пределе  в интегралах, а также  равномерных оценок на  $\partial \overline D(r)$ через $N_r$ на производные по радиусу функций $v_n$, переходя к пределу по $n\to +\infty$, из  
неравенства \eqref{VuMn} леммы \ref{lem2_2} получаем 
\begin{multline*}
\int_{\overline D(R)\setminus D(r)} V\dd \varDelta_u\leqslant 
\int_{\CC} V^*\dd \varDelta_u
\leqslant 
\int_{\CC} V^*\dd \varDelta_M+ \frac{r}{\pi}\int_0^{2\pi} \bigl(u(re^{i\theta})-M(re^{i\theta})\bigr)
\frac{\partial V}{\partial \vec n_{\operatorname{out}}}(re^{i\theta}) \dd \theta
\\
\overset{\eqref{V*}}{\leq} 
\int_{\overline D(R)\setminus D(r)} V\dd \varDelta_M+ \int_{\overline D(r)\setminus D(r^2/R)} V\biggl(\frac{r^2}{\bar z}\biggr)\dd \varDelta_M(z)+\frac{r}{\pi}\int_0^{2\pi} \bigl|u(re^{i\theta})-M(re^{i\theta})\bigr|N_r \dd \theta \\
\leq \int_{r<|z|\leq R} V\dd \varDelta_M+ \sup_{r^2/R\leq |z|\leq r} V\biggl(\frac{r^2}{\bar z}\biggr)
\varDelta_M \bigl(\overline D(r)\bigr)+\frac{r}{\pi}N_r \biggl(\int_0^{2\pi} \bigl|u(re^{i\theta})\bigr| \dd\theta +\int_0^{2\pi}\bigl|M(re^{i\theta})\bigr|\dd \theta\biggr)\\
\overset{\eqref{V*},\eqref{Cu}}{=}  \int_{r<|z|\leq R} V\dd \varDelta_M+ \sup_{r\leq |z|\leq R} V(z)\varDelta_M \bigl(\overline D(r)\bigr)
+\frac{r}{\pi}N_r\bigl(|u|^{\circ}(r)+|M|^{\circ}(r)\bigr),
\end{multline*}
что и даёт требуемую оценку \eqref{VuMnv}, завершая  доказательство  леммы \ref{lem2}. 
\end{proof}
Для применения леммы \ref{lem2} к убывающей последовательности функций \eqref{vn} с предельной функцией $V_R$ из \eqref{VR}
отметим, что функции $v_n$ удовлетворяют всем требованиям леммы \ref{lem2} по установленным выше перед леммой \ref{lem2_2} их свойствам и для них выполнены равномерные по $n$ неравенства для производных по радиусу
\begin{equation}\label{vnpr}
\biggl|\frac{\partial v_n}{\partial r} (re^{i\theta})\biggr|\overset{\eqref{vn}}{\leq} \sup_{\theta \in [0,2\pi]} k_n(\theta)
\Bigl|-\frac{1}{r^2}-\frac{1}{R^2}\Bigr|\leq \frac{2}{r^2} \sup_{\theta \in [0,2\pi]} k_1(\theta)=\frac{a}{r^2}=:N_r,
\end{equation}
где число $a:=\sup\limits_{\theta \in [0,2\pi]} k_1(\theta)$, очевидно, не зависит от $n\in \NN$. Таким образом, заключительная оценка 
\eqref{VuMnv} леммы \ref{lem2}  может быть записана для функций $V\overset{\eqref{VR}}{=}V_R$  и $M\overset{\eqref{uMs}}{\geqslant} u$ из \eqref{Ms} как 
\begin{multline*}
\int_{r<|te^{i\theta}|\leq R} \spf_{\bar S}(\theta)\Bigl(\frac{1}{t}-\frac{t}{R^2}\Bigr)\dd \varDelta_u(te^{i\theta})\overset{\eqref{VuMnv},\eqref{VR},\eqref{vnpr}}{\leqslant} 
\int_{r<|te^{i\theta}|\leq R} 
\spf_{\bar S}(\theta)\Bigl(\frac{1}{t}-\frac{t}{R^2}\Bigr)\dd \varDelta_M(re^{i\theta})
\\+\varDelta_M\bigl(\overline D(r)\bigr)\sup_{r\leq t\leq R}\Bigl(\frac{1}{t}-\frac{t}{R^2}\Bigr) \|\spf_{\bar S}\|_{\RR}
+ \frac{ar}{\pi r^2}\bigl(|u|^{\circ}(r)+|M|^{\circ}(r)\bigr).
\end{multline*}  
Отсюда, учитывая явный вид функции $M$ из \eqref{Ms} и её распределения масс Рисса из \eqref{rl}, имеем 
\begin{multline}\label{VuMnvRd}
\int_{r<|te^{i\theta}|\leq R} \frac{\spf_{\bar S}(\theta)}{t}\dd \varDelta_u(te^{i\theta})\overset{\eqref{Ms},\eqref{rl}}{\leqslant} \int_{r<|te^{i\theta}|\leq R} \frac{\spf_{\bar S}(\theta)t}{R^2}\dd \varDelta_u(te^{i\theta})+
\int_{r<|te^{i\theta}|\leq R} \frac{\spf_{\bar S}(\theta)}{t}\dd \varDelta_M(re^{i\theta})
\\+r\frac{l_{\conv \bar S}(2\pi)-l_{\conv \bar S}(0)}{2\pi} \frac{1}{r}\|\spf_{\bar S}\|_{\RR}
+ \frac{a}{\pi r}\bigl(|u|^{\circ}(r)+|M|^{\circ}(r)\bigr).
\end{multline}  
Первое, третье и четвёртое слагаемые из правой части этого неравенства оцениваются сверху числом, не зависящим от значений радиуса $r\geq r_0>0$. Действительно, для первого получаем  
\begin{equation*}\label{y}
\int_{r<|te^{i\theta}|\leq R} \frac{\spf_{\bar S}(\theta)t}{R^2}\dd \varDelta_u(te^{i\theta})\leq 
\|\spf_{\bar S}\|\frac{1}{R}\varDelta_u\bigl(\overline D(R)\bigr)\leq C_1,
\end{equation*}
где число $C_1\in \RR^+$ не зависит от $R\geq r_0$,  поскольку для субгармонической функции $u$ конечного типа при порядке $1$ выполнено \eqref{Du}. В третьем слагаемом $r$ просто исчезает  и оно оценивается сверху через некоторое  число $C_3\in \RR^+$.   
Наконец, 
$|M|^{\circ}(r)\overset{\eqref{Ms}}{\leq} \|\spf_{\bar S}\|r$ при всех $r\in \RR^+$, а для интегральных средних 
 $|u|^{\circ}(r)$ по окружностям $\partial \overline D(r)$ модуля субгармонической  функции 
$u\not\equiv -\infty$ конечного типа при порядке $1$ удовлетворяет, как следует, например, из \cite[лемма 6.2]{KhaShm19},
соотношению  $|u|^{\circ}(r)\underset{R\to +\infty}{=}O(r)$. Это даёт возможность оценить сверху четвёртое слагаемое  числом $C_4\in \RR^+$, не зависящим от $r\geq r_0>0$.
Таким образом, полагая $C:=C_1+C_3+C_4$ из \eqref{VuMnvRd} с учётом \eqref{areal} получаем оценку 
\begin{equation}\label{d}
\int_{r<|te^{i\theta}|\leq R} \frac{\spf_{\bar S}(\theta)}{t}\dd \varDelta_u(te^{i\theta})\overset{\eqref{VuMnvRd}}{\leqslant} 
\int_{r<|te^{i\theta}|\leq R} \frac{\spf_{\bar S}(\theta)}{t}\dd \varDelta_M(re^{i\theta})+C\quad \text{для всех $r\geq r_0$.}
\end{equation}  
 Это противоречит     равенству  \eqref{suMd}, что и завершает доказательство теоремы \ref{th1}.
\end{proof}

\begin{proof}[Доказательство следствия\/ {\rm \ref{cor1}}.]
Для ограниченной  односвязной области $S\subset \CC$ существует последовательность $(S_n)_{n\in \NN}$  компактов $S_n\subset S$
со связными дополнениями $\CC\setminus S_n$ при всех $n\in \NN$, объединение которых совпадает с  односвязной областью $S$. При этом для выпуклой оболочек $\conv S_n$ этих компактов  имеем $\area (\conv S_n)< \area (\conv S)$. Для полноты системы $\Exp^Z$ в $\Hol (S)$ с топологией равномерной сходимости на компактах  достаточно   показать, что система $\Exp^Z$ полна в каждом из пространств $C(S_n)\bigcap \Hol(\intr S_n)$  при  $n\in \NN$. Положим 
\begin{equation*}
d_n\underset{n\in \NN}{:=} \frac{1}{2}\bigl(\area (\conv S)- \area (\conv S_n)\bigr)>0.
\end{equation*}
При фиксированном $n\in \NN$ из равенства \eqref{sZar} следует существование возрастающей неограниченной последовательности $(a_k)_{k\in \NN}$ чисел $a_k>1$, для которой 
\begin{equation*}
\limsup_{0<r\to +\infty} \int_{r}^{a_kr}\frac{Z_{\spf_{\bar S}}^{\rad}(t)}{t^2}\dd t \geq \frac{1}{\pi} \bigl(\area(\conv  S_n)+d_n\bigr)\ln a_k \quad\text{при всех $k\in \NN$}.
 \end{equation*}
 Отсюда для каждого $k\in \NN$ найдётся достаточно большое  $r_k\geq 1$, для которого
\begin{equation*}
 \int_{r_k}^{a_kr_k}\frac{Z_{\spf_{\bar S}}^{\rad}(t)}{t^2}\dd t \geq \frac{1}{\pi} \bigl(\area(\conv  S_n)+d_n)\ln a_k-1
=\frac{1}{\pi} \area(\conv  S_n)\ln \frac{a_kr_k}{r_k}+\frac{1}{\pi} d_n\ln a_k-1,
\end{equation*}
что может быть записано как неравенства 
\begin{equation*}
 \int_{r_k}^{a_kr_k}\frac{Z_{\spf_{\bar S}}^{\rad}(t)}{t^2}\dd t- \frac{\area(\conv  S_n)}{\pi} \ln \frac{a_kr_k}{r_k}\geq\frac{1}{\pi} d_n\ln a_k-1,
\end{equation*}
Применяя здесь  операцию $\sup$ по $k$ к обеим частям, получаем 
 \begin{equation*}
 \sup_{k\in \NN}\biggl(\int_{r_k}^{a_kr_k}\frac{Z_{\spf_{\bar S}}^{\rad}(t)}{t^2}\dd t- \frac{\area(\conv  S_n)}{\pi} \ln \frac{a_kr_k}{r_k}\biggr)\geq \frac{d_n}{\pi}  \sup_{k\in \NN} \ln a_k-1=+\infty,
\end{equation*}
поскольку $d_n>0$. Тем более, имеет место равенство 
\begin{equation*}
\sup_{1<r<R<+\infty}\biggl(\int_{r}^{R}\frac{Z_{\spf_{\bar S}}^{\rad}(t)}{t^2}\dd t- 
\frac{\area(\conv  S_n)}{\pi} \ln \frac{R}{r}\biggr)=+\infty.
\end{equation*}
Отсюда по теореме \ref{th1} система $\Exp^{Z}$ полна в пространстве $C(S_n)\bigcap \Hol(\intr S_n)$. 
В силу произвола в выборе $n\in \NN$ получаем и полноту  системы $\Exp^{Z}$ в пространстве $\Hol(S)$, что и требовалось. 
\end{proof}

\end{document}